\input ./style/arxiv-general.cfg
\documentclass[aop,MSNbibl,seceqn,dvips]{arximspdf}
\makeatletter
   \@ifpackageloaded{graphicx}{}{\usepackage{graphicx}}
\makeatother

\doi{10.1214/15-AOP1008}
\volume{45}
\issue{1}
\pubyear{2017}
\firstpage{210}
\lastpage{224}
\docsubty{FLA}

\makeatletter
\newcommand{\lleft}{\left}
\newcommand{\rrvert}{\vert}
\newcommand{\rright}{\right}
\newcommand{\llvert}{\vert}

\newtheorem{Theorem}{Theorem}[section]
\newproclaim{Definition}[Theorem]{Definition}
\newtheorem{Proposition}[Theorem]{Proposition}
\newtheorem{Lemma}[Theorem]{Lemma}
\newproclaim{Remark}[Theorem]{Remark}
\makeatother

\begin{document}
\begin{frontmatter}

\title{Permanental vectors with nonsymmetric kernels}
\runtitle{Nonsymmetric permanental kernels}

\begin{aug}
\author[A]{\fnms{Nathalie}~\snm{Eisenbaum}\corref{}\ead
[label=e1]{nathalie.eisenbaum@upmc.fr}}
\runauthor{N. Eisenbaum}
\affiliation{CNRS, Universit\'e Pierre et Marie Curie}
\address[A]{Laboratoire de probabilit\'es\\\quad  et mod\`eles al\'eatoires\\
Universit\'e Paris
6\\
Case 188, 4 Place Jussieu\\ 75252 Paris cedex 05\\ France\\\printead{e1}}
\end{aug}

%
\received{\smonth{6} \syear{2014}}
%
\revised{\smonth{1} \syear{2015}}

%
\begin{abstract}
A permanental vector with a symmetric kernel and index $2$ is a squared
Gaussian vector.
The definition of permanental vectors is a natural extension of the
definition of squared
Gaussian vectors to nonsymmetric kernels and to positive indexes. The
only known permanental
vectors either have a positive definite kernel
or are infinitely divisible. Are there some others? We present a
partial answer to this question.
\end{abstract}

%
\end{frontmatter}

\section{Introduction}

A real-valued positive vector $\psi= (\psi_i, 1 \leq i \leq n)$
is a permanental vector if its Laplace transform satisfies
for every
$(\alpha_1, \alpha_2, \ldots,\alpha_n)$ in $\mathbb{R}^n_+$
%
\begin{equation}
\label{perm} \mathbb{E} \Biggl[\exp \Biggl\{-{1\over2} \sum
_{i = 1}^n \alpha_i \psi
_{i} \Biggr\} \Biggr] =\llvert I + \alpha G \rrvert ^{-\beta},
\end{equation}
where $I$ is the $n\times n$-identity matrix, $\alpha$ is the diagonal
matrix $\operatorname{Diag}((\alpha_i)_{1 \leq i \leq n})$,
$G = (G(i,j))_{1\leq i,j \leq n}$ and $\beta$ is a fixed positive number.

Such a vector $(\psi_i, 1 \leq i \leq n)$ is a permanental vector
with kernel
$(G(i,j),\allowbreak  1 \leq i,j \leq n)$ and index $\beta$. Note that the kernel
of $\psi$ is not uniquely determined. Indeed any matrix $D G D^{-1}$
with $D$  $n\times n$-diagonal matrix with nonzero entries is a
kernel for $\psi$. The matrices $G$ and $DGD^{-1}$ are said to be
\textit{diagonally equivalent}. But remark that $\psi$ also admits $G^t$ for
kernel. More generally, the kernels of $\psi$ are said to be \textit{effectively equivalent}.

Vere-Jones has established a necessary and sufficient condition on
the couple $(G, \beta)$ for the existence of such a vector. His
criterion is reminded at the beginning of Section~\ref{2}.

For $G$ $n\times n$-symmetric positive definite matrix and $\beta=
2$, (\ref{perm}) is the Laplace transform of the vector $(\eta_1^2,
\eta_2^2, \ldots, \eta_n^2)$ where $(\eta_1, \eta_2, \ldots, \eta_n)$ is
a centered Gaussian vector with covariance $G$. The definition of
permanental vectors hence represents an extension of the definition of
squared Gaussian vectors. The question is: to which point? More
precisely, one already knows two classes of matrices that satisfy
Vere-Jones criterion: the symmetric positive definite matrices and the
inverse $M$-matrices (a nonsingular matrix $A$ is a $M$-matrix if its
off-diagonal entries are nonpositive and the entries of $A^{-1}$ are
nonnegative). Up to effective equivalence, these are the only known
examples of permanental kernels. The question becomes:
Is there an irreducible permanental kernel that would not belong to any of this two classes?

\noindent This two classes correspond respectively to vectors with the Laplace
transform of a squared Gaussian vector to a positive power and to
infinitely divisible permanental vectors. Infinitely divisible
permanental processes are connected to local times of Markov processes
thanks to Dynkin's isomorphism theorem and its extensions (see \cite
{EK2}). Besides, we have shown in \cite{E} that for permanental
vectors, infinite divisibility and positive correlation are equivalent
properties.

\noindent In dimension one, obviously, the above two classes are identical and
the answer is negative. One easily checks that a 2-dimensional
permanental vector with index $2$ is a squared Gaussian couple.
Moreover, Vere-Jones \cite{V1}, solving a question raised by L\'evy \cite{L}, proved that a squared Gaussian couple is always
infinitely divisible. Hence, in this case also the two classes are
identical and the answer is negative. In dimension 3, the situation is
different. Indeed, Kogan and Marcus \cite{KM} have shown that if the
kernel of a $3$-dimensional permanental vector is not effectively
equivalent to a symmetric matrix (in short, is not
symmetrizable), then it is diagonally equivalent to an inverse
$M$-matrix. Since there exist inverse $M$-matrices that are not
symmetrizable, the two classes are not identical and have a nonempty
intersection. But the answer to the above question remains negative.

\noindent The case of dimension $d$ strictly greater than $3$ is still an open
question. We show here that if the kernel of a $d$-dimensional
permanental vector is strongly not symmetrizable, meaning that none of
its principal submatrices of dimension $3$ is symmetrizable, then it
is diagonally equivalent to the inverse of an $M$-matrix. The result
presented below is actually a little stronger and suggests that the
answer should still be negative in the general case. In other words,
one might think that the permanental vectors with a kernel not
effectively equivalent to a symmetric matrix, are always infinitely divisible.

\noindent For a set of indexes $I$, we adopt the notation: $G_{I \times I} =
(G(i,j))_{ (i,j) \in I \times I}$.

\begin{Theorem}\label{T1}
For $d > 3$, let $\psi$ be a $d$-dimensional
permanental vector with kernel $G$. Assume that there exists at most
one subset $I$ of three indexes such that $G_{I \times I} $ is
symmetrizable. Then $\psi$ is infinitely divisible.
\end{Theorem}

\noindent Theorem~\ref{T1} can also be stated as follows:

\vspace*{12pt}
\noindent\textit{For $\psi= (\psi_i)_{1 \leq i \leq d}$ permanental vector of
dimension $d>3$, assume that there exists at most three integers $1
\leq i_1, i_2,i_3 \leq d$ such that $(\psi_{i_1}, \psi_{i_2}, \psi
_{i_3})$ has the Laplace transform of a squared Gaussian vector to some
power. Then $\psi$ is infinitely divisible.}\vspace*{12pt}

\noindent The proof of Theorem~\ref{T1} is given in Section~\ref{2}. Section~\ref{1} introduces the needed preliminaries and definitions. Section~\ref{3} presents some examples and remarks.

\section{Preliminaries}\label{1}

We remind first the necessary and sufficient condition established
by Vere-Jones \cite{V}
for a given matrix $K$ to be the kernel of a permanental vector.

\textit{A $n\times n$-matrix $K$ is the kernel of a permanental vector
with index $\beta> 0$
if and only if}:
\begin{longlist}[(II)]
\item[(I)]\textit{ All the real eigenvalues of $K$ are nonnegative.}

\item[(II)] \textit{For every $\gamma> 0$, set $K_{\gamma} = (I + \gamma K)^{-1}
K$, then $K_{\gamma}$ is $\beta$-positive definite.}
\end{longlist}

A $n\times n$-matrix $M = (M(i,j))_{1\leq i,j \leq n}$ is said to be
$\beta$-positive definite if
for every integer $m$, every (not necessarily distinct) $k_1, k_2,\ldots,
k_m$ in $\{1,2,\ldots, n\}$
\[
\mathrm{per}_{\beta} \bigl( \bigl(M(k_i,k_j)
\bigr)_{1\leq i,j \leq m} \bigr) \geq0,
\]
where for any $m\times m$-matrix $A = (A(i,j))_{1 \leq i,j \leq m}$,
the quantity $\mathrm{per}_{\beta}(A)$ is defined as follows: $\mathrm{per}_{\beta}(A) = \sum_{\tau\in{\mathcal S}_m} \beta^{\nu(\tau)}
\prod_{i = 1}^m A_{i,\tau(i)}$,
with ${\mathcal S}_m$ the set of the permutations on $\{1,2,\ldots,m\}$, and
$\nu(\tau)$ the signature of $\tau$.

Note that the property of $\beta$-positive definiteness for a
matrix $M$ is supported by an infinite family of matrices derived from $M$.

\smallskip

The proposition below is just the regrouping of results of Kogan and
Marcus on the three dimensional permanental kernels. For the sake of
clarity, we explain where to find this results in \cite{KM}. Adopting
their convention, $0$ is both positive and negative.

\begin{Definition}\label{effequ} Two $d\times d$-matrices $A$ and $B$
are said to be effectively equivalent if for every $x$ in $\mathbb{R}^d\dvtx \llvert   I + x A\rrvert    =\llvert   I + x B\rrvert   $.
\end{Definition}

\begin{Definition}\label{Dsym} A squared matrix is symmetrizable if it
is effectively equivalent to a symmetric matrix.
\end{Definition}

\begin{Proposition}\label{P1} Let $\psi$ be a 3-dimensional
permanental vector with kernel $G = (G_{ij})_{1 \leq i,j \leq3}$. Then
we have:
\begin{longlist}[(iii)]
\item[(i)] $G$ is diagonally equivalent either to a matrix with all
positive entries or to a matrix with all negative off-diagonal entries.

\item[(ii)] If $G$ has all its off-diagonal entries strictly negative,
then $G$ is diagonally equivalent to a symmetric matrix.

\item[(iii)] If $G$ has all its off-diagonal entries strictly positive,
then $G$ is either an inverse $M$-matrix or it is diagonally equivalent
to a symmetric matrix.

\item[(iv)] If $G$ has one or more zero off-diagonal entries, it is
effectively equivalent to a symmetric matrix $\tilde{G}$ such that
$\tilde{G}_{ij} = 0$ when $G_{ij}G_{ji} = 0$.
\end{longlist}
\end{Proposition}

Up to some misprints, (i) is Remark~2.1 in \cite{KM}, which is a
consequence of the fact that $G_{ij} G_{ji} \geq0$ for every $1 \leq
i,j \leq3$. Indeed, for example, for $G$ with only positive entries
\[
\pmatrix{ G_{11} & G_{12} & G_{13}
\cr
G_{21} & G_{22}& G_{23}
\cr
G_{31} &
G_{32} & G_{33} }\quad \mbox{and}\quad \pmatrix{
G_{11} & -G_{12} & -G_{13}
\cr
-G_{21} &
G_{22}& G_{23}
\cr
-G_{31} & G_{32} &
G_{33} }
\]
are diagonally equivalent.
\begin{longlist}[(iii)]
\item[(ii)]  is established in the first part of the proof of Lemma~4.1 in
\cite{KM}.

\item[(iii)]  is established in the first part of the proof of Lemma~5.1 in
\cite{KM}.

\item[(iv)]  When $G$ is diagonally equivalent to a matrix with all negative
off-diagonal entries, this is a consequence of the second part of the
proof of Lemma~4.1 in \cite{KM}.
When $G$ is diagonally equivalent to a matrix with all positive
off-diagonal entries, (iv) is a consequence of the last paragraph of
the proof of Lemma~5.1 in \cite{KM}
together with its Lemma~2.3 cleaned from a misprint. For the last
sentence of Lemma~2.3 to be correct, the word ``diagonally'' should be
replaced by ``effectively.'' Indeed Lemma~2.3 in \cite{KM}, assuming
that the two matrices
\[
\pmatrix{ 1 & 0 & c_2
\cr
a_2 & 1& b_1
\cr
c_1 & 0 & 1 }\quad \mbox{and}\quad %
\pmatrix{ 1 & 0 &
\sqrt{c_1c_2}
\cr
0 &1& 0
\cr
\sqrt{c_1c_2}
& 0 & 1 }
\]
are permanental kernels, states that for $a_2b_1 c_1c_2 \neq 0$, they
are diagonally equivalent. But they cannot be diagonally equivalent.
They are effectively equivalent.
\end{longlist}


We will use repeatedly the following lemma which is an elementary remark.

\begin{Lemma}\label{L1}
For $A$ $n\times n$-matrix, the following
points are equivalent:
\begin{longlist}[(iii)]
\item[(i)] $A$ is diagonally equivalent to a symmetric matrix.

\item[(ii)] For every couple $(D_1, D_2)$ of diagonal $n\times n$-matrices
with strictly positive diagonal entries, $D_1AD_2$ is diagonally
equivalent to a symmetric matrix.

\item[(iii)] There exist two diagonal $n\times n$-matrices with strictly
positive diagonal entries $D_1$ and $D_2$ such that $D_1 A D_2$ is
diagonally equivalent to a symmetric matrix.
\end{longlist}
\end{Lemma}

\begin{Remark}\label{R0}
 In view of Proposition~\ref{P1}, a
permanental kernel of dimension~$3$, $G = (G_{ij})_{1 \leq i,j \leq
3}$, is symmetrizable iff:
\begin{itemize}
\item either $G$ has an off-diagonal entry equal to zero,
\item either $G$ has no zero off-diagonal entry and it is diagonally
equivalent to a symmetric matrix with strictly positive entries.
\end{itemize}
To check whether a $3\times3$-matrix $K$ without zero off-diagonal
entry, is symmetrizable, one has first to check the existence of a
signature matrix $\sigma$ (a diagonal matrix with $\llvert  \sigma_{ii}\rrvert    = 1,
1\leq i \leq3$) such that:
$\sigma K \sigma= (\llvert  K_{ij}\rrvert   )_{1 \leq i,j \leq3}$, and then check that
\[
\bigl\llvert K(1,2) K(2,3) K(3,1) \bigr\rrvert = \bigl\llvert K(2,1) K(1,3)
K(3,1) \bigr\rrvert.
\]
\end{Remark}

\section{Proof of Theorem \texorpdfstring{\protect\ref{T1}}{1.1}}
\label{2}

\textit{Step~1}: Assume that $d = 4$ and that $G$ has no symmetrizable
$3 \times3$-principal submatrices, we show then that $\psi$ is
infinitely divisible.

Thanks to Remark~\ref{R0}, we know that
$G$ has no entry equal to $0$. Moreover, in view of (i) and (ii) of
Proposition~\ref{P1}, every $3\times3$-principal submatrix of $G$
has to be diagonally equivalent to a matrix with all positive
entries. This means that for every subset of three indexes $I$,
there exists $S_I$ from $I$ into $\{-1, +1\}$ such that $S_I(i) G(i,j)
S_I(j) \geq0$, for $i,j \in I$.
This leads to
\[
G(i,j)G(j,k)G(k,i) > 0\qquad \forall i,j,k \in\{1,2,3,4\}.
\]
Since $G$ has no zero entry, this property implies the existence of $S$
from $\{1,2,3,4\}$ into $\{-1, +1\}$ such that: $S(i) G(i,j) S(j) > 0,
\forall i,j \in\{1,2,3,4\}$.
Since $\psi$ also admits for kernel $SGS$, we will
assume from now that the entries of $G$ are all strictly positive.

For $\sigma> 0$,
consider the $3$-dimensional vector $\phi_{\sigma}$ with Laplace transform
%
\begin{equation}
\label{LT} {\mathbb{E}[ \exp\{-{(1/2)} \sum_{j = 1}^{3}\lambda_j \psi_j\}
\exp\{
-{(\sigma/2)} \psi_4\}] \over\mathbb{E}[ \exp\{-{(\sigma/
2)} \psi
_4\}]}.
\end{equation}
This vector is a permanental vector with the same index as $\psi$ and
admits for kernel
$H({\sigma}, G)$ (see \cite{KM})
\[
H({\sigma}, G) = \biggl(G(i,j) - {\sigma\over1 + \sigma G(4,4)} G(i,4) G(4,j)
\biggr)_{ 1 \leq i,j \leq3}.
\]
For which values of $\sigma$, is $H(\sigma, G)$ symmetrizable?
For $\sigma> 0$, we have\break ${\sigma\over1 + \sigma G(4,4)} < {1
\over G(4,4)}$.
We set: $\Gamma= ({G(i,j)\over G(i,4)G(4,j)})_{1 \leq i,j \leq4}$.
Making use of Lemma~\ref{L1}, we are looking for the values of $c $ in
$(0, {1 \over G(4,4)})$ such that $(\Gamma(i,j) - c)_{1 \leq i,j \leq
3}$ is symmetrizable. In view of Remark~\ref{R0}, this can occur in
two ways:
\begin{itemize}
\item either $(\Gamma(i,j) - c)_{1 \leq i,j \leq3}$ has an off-diagonal
entry equal to zero.

\item either $(\Gamma(i,j) - c)_{1 \leq i,j \leq3}$ has no zero
off-diagonal entry and it is diagonally equivalent to a symmetric matrix.
\end{itemize}

\noindent The first possibility is excluded because it would imply the
existence of $i$ and $j$ distinct from $4$, such that: ${G(i,j)\over
G(i,4) G(4,j)} < {1 \over G(4,4)}$. But since the entries of $G$ are
all strictly positive, we know by assumption that $G_{\{i,j,4\} \times
\{i,j,4\}}$ is an inverse $M$-matrix. This last property implies in
particular that: $G(i,j) G(4,4) \geq G(i,4)G(4,j)$. This can bee seen
by computing the inverse of $G_{\{i,j,4\} \times\{i,j,4\}}$ or by
using Willoughby's paper \cite{W}.

We now study the second possibility. Since $(\Gamma(i,j)
-c)_{1\leq i,j \leq3}$ has only strictly positive entries, we know,
thanks to Lemma~\ref{L1}(iii), that $(\Gamma(i,j) - c)_{1 \leq i,j
\leq3}$ \mbox{is diagonally} equivalent to a symmetric matrix if and only
if\break
$({\Gamma(i,j) - c \over(\Gamma(i,3) -c)(\Gamma(3,j) - c)})_{ 1
\leq i,j \leq3}$ is. Denote this last matrix by $A_c$. Since $A_c(i,3)
= A_c(3,j) = {1\over\Gamma(3,3) - c}$ for every $1\leq i,j \leq3$,
$A_c$ is diagonally equivalent to a symmetric matrix if and only if
$(A_c(i,j))_{1 \leq i,j \leq2}$ is symmetric.
This translates into
\[
{\Gamma(1,2) - c \over(\Gamma(1,3) -c)(\Gamma(3,2) - c)} = {\Gamma
(2,1) - c \over(\Gamma(2,3) -c)(\Gamma(3,1) - c)},
\]
which means that $c$ must solve a polynomial equation with degree $3$.
Hence, only the
two following cases might occur:
\begin{itemize}
\item either there are at most three distinct values for $c$ such that
$(\Gamma(i,j) - c)_{1 \leq i,j \leq3}$ is diagonally equivalent to a
symmetric matrix,

\item either for every value of $c$, $(\Gamma(i,j) - c)_{1 \leq i,j
\leq3}$ is diagonally equivalent to a symmetric matrix.
\end{itemize}
In the later case, one obtains in particular $(\Gamma(i,j))_{1 \leq
i,j \leq3}$ is diagonally equivalent to a symmetric matrix. Thanks to
Lemma~\ref{L1}, this implies that $(G(i,j))_{1 \leq i,j \leq3}$ is
diagonally equivalent to a symmetric matrix. But this is excluded by assumption.

\noindent Consequently, except for at most three distinct values of $\sigma$,
$H(\sigma, G)$ is not symmetrizable.

Set now $G_{\sigma} = (I + \sigma G)^{-1} G$. We have shown
(Proposition~3.2 in \cite{E}) that there exists a permanental vector
$\psi_{\sigma}$ with the same index as $\psi$, admitting $G_{\sigma
}$ for kernel and such that its Laplace transform satisfies
\begin{eqnarray*}
&&\mathbb{E} \Biggl[ \exp\Biggl\{-{1\over2} \sum
_{j = 1}^{4}\lambda_j \psi
_{\sigma
}(j) \Biggr\} \Biggr]
\\
&&\qquad = \mathbb{E} \Biggl[{ \exp \{- {(\sigma/2)} \sum_{i = 1}^4 \psi
(i) \}
\over\mathbb{E}[\exp\{- {(\sigma/2)} \sum_{i = 1}^4 \psi(i) \}
]} \exp \Biggl\{-
{1\over2} \sum_{j = 1}^{4}
\lambda_j \psi(j) \Biggr\} \Biggr].
\end{eqnarray*}
Hence, thanks to (\ref{LT}), it also satisfies
%
\begin{equation}
\label{LT2} \mathbb{E} \Biggl[ \exp \Biggl\{-{1\over2} \sum
_{j = 1}^{3}\lambda_j \psi
_{\sigma
}(j) \Biggr\} \Biggr] = c(\sigma) \mathbb{E} \Biggl[ \exp \Biggl\{-
{1\over2} \sum_{j = 1}^{3}(
\lambda_j + \sigma) \phi _{\sigma}(j) \Biggr\} \Biggr],
\end{equation}
where $c(\sigma) = { \mathbb{E}[ \exp\{- {(\sigma/2)} \psi(4)\}
]
\over\mathbb{E}[ \exp\{- {(\sigma/2)} \sum_{i = 1}^4 \psi(i)\}
]} $.

We show now that $(\psi_{\sigma}(i))_{1 \leq i \leq3}$ admits a
kernel without zero entry.
To do so, we first adopt the notation: ${\underline G}_{\sigma} =
(G_{\sigma}(i,j))_{1 \leq i,j \leq3}$.
Note that (\ref{LT2}) can be written as follows for every $x \in
\mathbb{R}^3$:
\[
\llvert I + x {\underline G}_{\sigma}\rrvert = {\llvert  I + ( x + \sigma I) H_{\sigma}\rrvert
\over\llvert  I + \sigma H_{\sigma}\rrvert   }.
\]
Besides denote by $R_{\alpha}$ the $\alpha$-resolvent of $H_{\sigma
}\dvtx R_{\alpha} = (I + \alpha H_{\sigma})^{-1} H_{\sigma}$, then we have
\begin{eqnarray*}
\llvert  I + x R_{\sigma}\rrvert    &=&\llvert  I + x (I + \sigma H_{\sigma})^{-1} H_{\sigma}\rrvert    =\llvert   I + \sigma H_{\sigma}\rrvert   ^{-1}\llvert   I + (x + \sigma I) H_{\sigma}\rrvert\\
&=& \llvert  I +x {\underline G}_{\sigma}\rrvert .
\end{eqnarray*}

\noindent Consequently: ${\underline G}_{\sigma} $ and $R_{\sigma}$ are
effectively equivalent. This implies that ${(\underline G}_{\sigma
})^{-1}$ and $R_{\sigma}^{-1}$ are effectively equivalent.

Assume that ${\underline G}_{\sigma}$ has an off-diagonal entry
equal to zero. Then thanks to Proposition~\ref{P1}(iv), ${\underline
G}_{\sigma}$ is effectively equivalent to a symmetric matrix. This
implies that ${(\underline G}_{\sigma})^{-1}$, and consequently
$R_{\sigma}^{-1}$, is effectively equivalent to a symmetric matrix. Since
$R_{\sigma}^{-1} = H_{\sigma}^{-1} + \sigma I$, we obtain that
$(H_{\sigma}^{-1} + \sigma I)$ is effectively equivalent to a
symmetric matrix, which easily implies that $H_{\sigma}^{-1}$, and
then $H_{\sigma}$ must be effectively equivalent to a symmetric
matrix. Except for at most three values of $\sigma$, this is not true.

Consequently, we have obtained, except for at most three values of
$\sigma$, that $(G_{\sigma}(i,j))_{1 \leq i,j \leq3}$ has no zero
entry and is not effectively equivalent to a symmetric matrix. In view
of Proposition~\ref{P1}, the only possibility for $(G_{\sigma
}(i,j))_{1 \leq i,j \leq3}$ is to be diagonally equivalent to an
inverse $M$-matrix.

Hence, there exists a function $S$ from $\{1,2,3\}$ into $\{-1, +1\}
$ such that for every $i,j$ in $\{1,2,3\}$:
\[
S(i) G_{\sigma}(i,j) S(j) > 0,
\]
which leads to
%
\begin{equation}
\label{cle} G_{\sigma}(i,j) G_{\sigma}(j,k) G_{\sigma}(k,i)
> 0
\end{equation}
for every $i,j,k$ in $\{1, 2, 3 \}$.

The choice of the three indexes $1$, $2$ and $3$, being arbitrary,
we conclude that excepted for at most a finite number of values of
$\sigma$, $G_{\sigma}$ has no zero entry
and satisfies
\[
G_{\sigma}(i,j) G_{\sigma}(j,k) G_{\sigma}(k,i) > 0
\]
for every $i,j,k$ in $\{1, 2, 3, 4 \}$.

This last property implies that there exists a function ${S}_{\sigma
}$ from $\{1,2,3,4\}$ into $\{-1, +1\}$ such that for every $i$, $j$ in
$\{1,2,3,4\}$
%
\begin{equation}
\label{signe} {S}_{\sigma}(i) G_{\sigma}(i,j) {S}_{\sigma}(j)
> 0.
\end{equation}

\noindent For the three values of $\sigma$ that we have excluded, we still
have (\ref{signe}). Indeed
assume that there exists such a value and that it is strictly positive.
Denote it by~$\alpha$.
We know now that there exists $\varepsilon> 0$, such that for every
$\sigma$ in $(\alpha, \alpha+ \varepsilon]$,
$G_{\sigma}$ satisfies (\ref{signe}). Note that $G_{\sigma}$ is the
$(\sigma- \alpha)$-resolvent of $G_{\alpha}$:
\[
G_{\alpha} = G_{\sigma} \bigl( I - (\sigma- \alpha)G_{\sigma}
\bigr)^{-1}.
\]
For $(\sigma- \alpha)$ small enough, we have
$G_{\alpha} = \sum_{k = 1}^{\infty}(\sigma- \alpha)^kG_{\sigma
}^k$. Making use of (\ref{signe}), one obtains for every $1\leq i,j
\leq4$:
\[
S_{\sigma}G_{\alpha}S_{\sigma}(i,j) = \sum
_{k = 1}^{\infty}(\sigma- \alpha)^kS_{\sigma}G_{\sigma}^kS_{\sigma}(i,j)
= \sum_{k = 1}^{\infty}(\sigma-
\alpha)^k(S_{\sigma}G_{\sigma}S_{\sigma
})^k
(i,j) > 0.
\]

\noindent Obviously, (\ref{signe}) implies that $G_{\sigma}$ is $\beta
$-positive definite for every $\beta> 0$.
Vere-Jones criteria allows then to conclude that $\psi$ is infinitely
divisible.

\noindent \textit{Conclusion of Step~1}:
We have actually established that if $G$ is a permanental kernel of
dimension $4$, with no
symmetrizable $3 \times3$-principal submatrices, then it is the
kernel of an infinitely divisible permanental vector, and moreover,
for every $\sigma> 0$, its $\sigma$-resolvent $G_{\sigma}$ has a no
zero entry.

\smallskip

\textit{Step~2}:   Define the claim $(R_n)$ as follows.

$(R_n)$: If $G$ is a $n\times n$-square matrix without
symmetrizable $3\times3$-principal submatrix, then $\psi$ is
infinitely divisible and for every $\sigma>0$, its $\sigma$-resolvent
$G_{\sigma}$ has no zero entry.

We have just established $(R_4)$.
Assume that $(R_n)$ is satisfied. We now
establish $(R_{ n+1})$.
First note that $G$ is diagonally equivalent to a matrix with only
strictly positive entries. Indeed, using exactly the same argument as
at the beginning of Step~1, one obtains
\[
G(i,j)G(j,k)G(k,i) > 0\qquad \forall i,j,k \in\{1,2,\ldots,n\}.
\]
Similarly, as in Step~1, one concludes that there exists $S$ from $\{
1,2,\ldots,n\}$ into $\{-1,+1\}$ such that $S(i) G(i,j) S(j) > 0, \forall
i,j \in\{1,2,\ldots,n\}$.

Hence, we can assume that all the entries of $G$ are strictly
positive. Then
consider the vector $(\phi_{\sigma}(i))_{1\leq i \leq n}$ with
Laplace transform
%
\begin{equation}
\label{LTH} {\mathbb{E}[ \exp\{-{(1/2)} \sum_{j = 1}^{n}\lambda_j \psi_j\}
\exp\{
-{(\sigma/2)} \psi_{n+1}\}] \over\mathbb{E}[ \exp\{-{(\sigma
/2)}
\psi_{n+1}\}]}.
\end{equation}
The vector $(\phi_{\sigma}(i))_{1\leq i \leq n}$ is a permanental
vector admitting for kernel $H(\sigma, G)$ defined by
\[
H({\sigma}, G) = \biggl(G(i,j) - {\sigma\over1 + \sigma G(n+1,n+1)} G(i,n+1) G(n+1,j)
\biggr)_{ 1 \leq i,j \leq n}.
\]
We look for the values of $\sigma$ such that $H(\sigma, G)$ would
have a $3\times3$-principal symmetrizable matrix. We set $\Gamma= (
{G(i,j)\over G(i,n+1) G(n+1,j)})_{1 \leq i,j \leq n+1}$. We hence look
for the values of $c$ in $(0, {1\over G(n+1,n+1)})$ such that $(\Gamma
(i,j) - c)_{1 \leq i,j \leq n}$ would have a symmetrizable $3\times
3$-principal submatrix. We fix $I$, a subset of three elements of $\{
1,2,\ldots,n\}$.
Similarly, as in the case $d = 4$, we know that $(\Gamma- c)_{I\times
I}$ has no zero entry. The only way for $(\Gamma- c)_{I\times I}$ to
be symmetrizable is to be diagonally equivalent to a symmetric matrix.
Again as in Step~1, we have:
\begin{itemize}
\item either $(\Gamma- c)_{I\times I}$ is symmetrizable for at most
three distinct values of $c$,

\item either for every real $c$, $(\Gamma- c)_{I\times I}$
is diagonally equivalent to a symmetric matrix.
\end{itemize}
In the second case, one obtains that $\Gamma_{I\times I}$, and
consequently $G_{I\times I}$ thanks to Lemma~\ref{L1}, is diagonally
equivalent to a symmetric matrix, which is excluded by assumption.
Consequently, for $\sigma$ outside of a finite set, $H(\sigma, G)$
does not contain any $3\times3$-principal symmetrizable matrix.
Thanks to our assumption on the case $d= n$, we know that $\phi
_{\sigma}$ is infinitely divisible and that for every $\alpha>0$,
$R_{\alpha}$, the $\alpha$-resolvant of $H_{\sigma}$ has no zero
entry. Besides there exists a permanental vector with the same index as
$\phi_{\sigma}$, admitting $R_{\alpha}$ for kernel (see Proposition~3.2 in \cite{E}). Making use of Vere-Jones criteria, one easily shows
that this permanental vector is infinitely divisible too.

Setting $G_{\sigma} = (I+ \sigma G)^{-1} G$, we know (Proposition~3.2 in \cite{E}) that there exists a permanental vector $\psi_{\sigma
}$ with the same index as $\psi$, admitting $G_{\sigma}$ for kernel
and such that its Laplace transform satisfies
\begin{eqnarray*}
&&\mathbb{E} \Biggl[ \exp \Biggl\{-{1\over2} \sum
_{j = 1}^{n+1}\lambda_j \psi
_{\sigma
}(j) \Biggr\} \Biggr]
\\
&&\qquad = \mathbb{E} \Biggl[{ \exp\{- {(\sigma/2)} \sum_{i = 1}^{n+1}
\psi(i) \}
\over\mathbb{E}[\exp\{- {(\sigma/2)} \sum_{i = 1}^{n+1} \psi
(i) \}]} \exp \Biggl\{-
{1\over2} \sum_{j = 1}^{n+1}
\lambda_j \psi(j) \Biggr\} \Biggr].
\end{eqnarray*}
Hence, thanks to (\ref{LTH}), it also satisfies
%
\begin{equation}
\mathbb{E} \Biggl[ \exp \Biggl\{-{1\over2} \sum
_{j = 1}^{n}\lambda_j \psi
_{\sigma
}(j) \Biggr\} \Biggr] = c(\sigma) \mathbb{E} \Biggl[ \exp \Biggl\{-
{1\over2} \sum_{j = 1}^{n}(
\lambda_j + \sigma) \phi _{\sigma}(j) \Biggr\} \Biggr],
\end{equation}
where $c(\sigma) = { \mathbb{E}[ \exp\{- {(\sigma/2)} \psi
(n+1)\}]
\over\mathbb{E}[ \exp\{- {(\sigma/2)} \sum_{i = 1}^{n+1} \psi
(i)\}]} $.

Similarly, as in Step~1, one shows that the two $n\times n$-matrices
$R_{\sigma}$ and
$(G_{\sigma}(i,j))_{1 \leq i,j \leq n}$ are effectively equivalent.
Note that for every $1 \leq i,j \leq n$:  $G_{\sigma}(i,j)G_{\sigma
}(j,i) = R_{\sigma}(i,j) R_{\sigma}(j,i)$.

For $\sigma$ outside a finite set, $R_{\sigma}$ has no zero entry, and hence
neither\break
$(G_{\sigma}(i,j))_{1 \leq i,j \leq n}$. Moreover, we know also
that $(\psi_{\sigma}(i))_{1 \leq i \leq n}$ is infinitely divisible.
In particular for every triplet of indexes $i$, $j$ and $k$ in $\{
1,2,\ldots,n\}$,
$(\psi_{\sigma}(i), \psi_{\sigma}(j), \psi_{\sigma}(k))$ is
infinitely divisible. Consequently,
$(G_{\sigma})_{\{i,j,k\}\times\{i,j,k\}}$ is diagonally equivalent to
an inverse $M$-matrix. The choice of the index $(n+1)$ being arbitrary,
we actually obtain that for $\sigma$ outside of a finite set there
exists a function ${S}_{\sigma}$ from $\{1,2,\ldots,n,n+1\}$ into $\{-1,
+1\}$ such that for every $i$, $j$ in $\{1,2,\ldots,n,n+1\}$
%
\begin{equation}
\label{signe1} {S}_{\sigma}(i) G_{\sigma}(i,j) {S}_{\sigma}(j)
> 0.
\end{equation}
For $\sigma$ element of the finite set of excluded values, one shows
that (\ref{signe1}) is still true exactly as we did it for (\ref
{signe}) in Step~1.

We conclude that for every $\sigma> 0$, $G_{\sigma}$ has no zero
entry and is $\beta$-positive definite for every $\beta> 0$. Thanks
to Vere-Jones criteria, $\psi$ is infinitely divisible and $(R_{n+1})$
is established.

\smallskip

\textit{Step~3}: Assume that $d= 4$ and that $G$ is such that the
matrices\break $(G(i,j))_{i,j \in\{1,2,4\}}$, $(G(i,j))_{i,j \in\{1,3,4\}}$
and $(G(i,j))_{i,j \in\{2,3,4\}}$ are not\vadjust{\goodbreak} symmetrizable. We show that
$\psi$ is infinitely divisible and that for every $\sigma> 0$, its
$\sigma$-resolvent $G_{\sigma}$ has no zero entry.

First, note that according Remark~\ref{R0}, the three matrices
$(G(i,j))_{i,j \in\{1,2,4\}}$, $(G(i,j))_{i,j \in\{1,3,4\}}$ and
$(G(i,j))_{i,j \in\{2,3,4\}}$ have no zero entry. Hence, $G$ has no
zero entry. Since these three matrices are all diagonally equivalent to
inverse of $M$-matrices, we can then easily establish the existence of
$S$ from $\{1,2,3,4\}$ into $\{-1, +1\}$ such that: $S(i) G(i,j) S(j) >
0, \forall i,j \in\{1,2,3,4\}$. We can hence assume that the entries
of $G$ are all strictly positive.

We now make the notation for $H(\sigma, G)$ more precise, by writing
\[
H( \sigma, G, 4) = \biggl(G(i,j) - {\sigma\over1 + \sigma
G(4,4)}G(i,4)G(4,j)
\biggr)_{1 \leq i,j \leq3}.
\]
Similarly, for any $k$ in $\{1,2,3,4\}$, $H(\sigma, G, k)$ is defined by
\[
H( \sigma, G, k) = \biggl(G(i,j) - {\sigma\over1 + \sigma
G(k,k)}G(i,k)G(k,j)
\biggr)_{ i,j \in\{1,2,3,4\}\setminus\{k\} }.
\]

Making use of the argument developed in Step~1, we know that for
each of the three matrices $H(\sigma, G, 3)$, $H(\sigma, G, 2)$ and
$H(\sigma, G, 1)$, there are at most three distinct values of $\sigma
$ for which they are not symmetrizable. Consequently, for $\sigma$
outside of a finite set, the three matrices $(G_{\sigma}(i,j))_{i,j
\in\{1,2,4\}}$,
$(G_{\sigma}(i,j))_{i,j \in\{1,3,4\}}$ and $(G_{\sigma}(i,j))_{i,j
\in\{2,3,4\}}$ have no zero entry and are diagonally equivalent to
inverse $M$-matrices. Setting $I_3 = \{1,2,4\}$, $I_2 = \{1,3,4\}$ and
$I_1 = \{2,3, 4\}$, we hence know that
there exist three functions $S_{3}, S_{2}$ and $S_{1}$ from
respectively $I_3$, $I_2$ and $I_1$ into $\{-1, +1\}$ such that for
every $p = 1,2$ or $3$, and every couple $(i,j)$ of $I_p$, we have
\[
S_p(i) G_{\sigma}(i,j) S_p(j) > 0.
\]
To determine the sign of $G_{\sigma}(1,2) G_{\sigma}(2,3) G_{\sigma
}(3,1)$, note that it has the same sign as
$ S_3(1)S_3(2)\cdot S_1(2)S_1(3)\cdot S_2(3)S_2(1)$. But
$S_3(1)S_2(1)$ has the same sign as $S_3(4)S_2(4)G_{\sigma}(4,1)^2$;
$S_3(2)S_1(2)$ has the same sign as $S_3(4)S_1(4)G_{\sigma}(2,4)^2$
and $S_1(3)S_2(3)$ has the same sign as $S_1(4)S_2(4)G_{\sigma
}(4,3)^2$. One obtains
\[
G_{\sigma}(1,2)\* G_{\sigma}(2,3) G_{\sigma}(3,1) > 0.
\]
 Consequently, for every $i$, $j$, $k$ in $\{1,2,3,4\}$ we have
\[
G_{\sigma}(i,j) G_{\sigma}(j,k) G_{\sigma}(k,i) > 0,
\]
which leads to the existence of a function $S$ from $\{1,2,3,4\}$ to $\{
-1,+1\} $ such that for every $i$,$j$ in $\{1,2,3,4\}$:
%
\begin{equation}
\label{signe2} S(i) G_{\sigma}(i,j) S(j) > 0.
\end{equation}
For $\sigma$ element of the finite set of excluded values, one shows
that (\ref{signe2}) is still true exactly as we did it for (\ref
{signe}) in Step~1.

We conclude that for every $\sigma>0$, $G_{\sigma}$ has no zero
entry and that $\psi$ is infinitely divisible.

\smallskip

\textit{Step~4}: We assume that $G$ has exactly one symmetrizable
principal $3\times3$-submatrix. Denote by $I$ the subset of the
corresponding three distinct indexes. We show that $\psi$ is
infinitely divisible and that for every $\sigma> 0$, its $\sigma
$-resolvent $G_{\sigma}$ has no zero entry.
Define the claim $(\tilde{R}_n)$ as follows.

$(\tilde{R}_n)$: If $G$ is a $n\times n$-square matrix with
exactly one symmetrizable $3\times3$-principal submatrix, then $\psi
$ is infinitely divisible and for every $\sigma> 0$, its $\sigma
$-resolvent $G_{\sigma}$ has no zero entry.

We just established $(\tilde{R}_4)$.
Assume that $(\tilde{R}_n)$ is satisfied we show now that $(\tilde
{R}_{n+1})$ is satisfied.

As in Step~3, one shows that we can assume that the entries of $G$
are strictly positive.
Note that for every index $p$ in $\{1,2,\ldots, n+1\}$, $H(\sigma, G, p)$
is the kernel of a $n$-dimensional permanental vector.
We still set $\Gamma=\break  ({G(i,j) \over G(i,p) G(p, j)})_{1 \leq i,j \leq
n+1}$. Fix $J$ subset of three elements of $\{1,2,\ldots, n+1\} \setminus
\{p\}$.
We look for the values $c$ in $(0, {1\over G(p,p)})$ for which $(\Gamma
(i,j) - c)_{i,j \in J\times J}$ is symmetrizable. Unless $J = I$, we
know, similarly as in Step~2, that $(\Gamma(i,j) - c)_{i,j \in J\times
J}$ has no off-diagonal zero entry. Hence, for $J \neq I$, the only
way for $(\Gamma(i,j) - c)_{i,j \in J\times J}$ to be symmetrizable is
to be diagonally equivalent to a symmetric matrix. We know that:
\begin{itemize}
\item either $(\Gamma(i,j) - c)_{i,j \in J\times J}$ is diagonally
equivalent to a symmetric matrix for at most three distinct values of
$c$,

\item either $(\Gamma(i,j) - c)_{i,j \in J\times J}$ is diagonally
equivalent for every value of $c$.
\end{itemize}
In the later case, one obtains $G_{J\times J}$ is symmetrizable, which
implies that $J = I$.
Consequently
for every $p$, and every $\sigma$ outside of a finite set, $H(\sigma,
G,p)$ contains at most one $3\times3$-symmetrizable principal
submatrix. If there is none, then Step~2 tells us that the
corresponding permanental vector is infinitely divisible and that for
every $\alpha> 0$, its $\alpha$-resolvent has no zero entry. If
$H(\sigma, G,p)$ has exactly one $3\times3$- symmetrizable principal
submatrix, we obtain the same property thanks to $(\tilde{R}_n)$.
Making use of the argument developed in Step~2, one shows that $\psi$
is infinitely divisible and for every $\sigma> 0$, its $\sigma
$-resolvent has no zero. We have hence
obtained $\tilde{R}_{n+1}$. This completes the proof of Theorem~\ref{T1}.

\section{Remarks and examples}\label{3}

%
\begin{Remark}
 Theorem~\ref{T1} can be reformulated in terms of
linear algebra as follows.

\textit{For $d > 3$, let $G$ be a $d\times d$-matrix with no zero
entry such that at most one of its $3\times3$-principal submatrices
is diagonally equivalent to a symmetric matrix. Assume that}:
\begin{longlist}[(II)]

\item[(I)] \textit{all the real eigenvalues of $G$ are nonnegative},

\item[(II)] \textit{there exists $\beta> 0$, such that for every $\gamma> 0$,
setting $G_{\gamma} = (I + \gamma G)^{-1} G$, $G_{\gamma}$ is $\beta
$-positive definite},
\end{longlist}
\textit{then $G$ is diagonally equivalent to an inverse $M$-matrix}.\vadjust{\goodbreak}

Assumptions (I) and (II) are necessary to obtain the conclusion. Indeed,
consider the following nonsingular matrix borrowed from \cite{JS}:
\[
A = \pmatrix{ 1 & 0,10 & 0,40 & 0,30
\cr
0,40 &1 & 0,40 & 0,65
\cr
0,10 & 0,20 &
1 & 0,60
\cr
0,15 & 0,30 & 0,60 & 1 }.
\]
It is not an inverse $M$-matrix, since $A^{-1}(2,3)$ is positive. But
note that every $3\times3$-principal submatrix is an inverse
$M$-matrix and is not symmetrizable. Consequently, $A$ is not the
kernel of a permanental vector.
\end{Remark}
%

\begin{Remark}
 The condition required by Theorem~\ref{T1} to
obtain infinite divisibility, is sufficient and not necessary. Indeed,
there exist nonsymmetrizable inverse $M$-matrices with more than one
symmetrizable $3\times3$-principal submatrix. Here is a family of
such matrices with dimension $4$:
\[
\Gamma= \pmatrix{ \Gamma(1,1) & a & a & \Gamma(4,4)
\cr
b &\Gamma(2,2) & e &
\Gamma(4,4)
\cr
b & e & \Gamma(3,3) & \Gamma(4,4)
\cr
\Gamma(4,4) & \Gamma(4,4) &
\Gamma(4,4) & \Gamma(4,4) },
\]
with $\Gamma(i,i) > e$ for $i = 1,2,3$;  $a,b, e > \Gamma(4,4)$ and
$e > a,b$.

\noindent For $a \neq b$, $\Gamma$ is not symmetrizable and has exactly two
symmetrizable $3\times3$-principal submatrices, $\Gamma_{\llvert  _{\{1,2,3\}\times\{1,2,3\}}}$ and $\Gamma_{\rrvert   _{\{2,3,4\}\times\{2,3,4\}}}$, and
two nonsymmetrizable $3\times3$-principal submatrices.

Here are now examples of matrices illustrating Theorem~\ref{T1}.
\[
\mbox{Set: } K = \pmatrix{ K(1,1) & e & a &K(4,4)
\cr
b & K(2,2) & a & K(4,4)
\cr
b & e & K(3,3)& K(4,4)
\cr
K(4,4)& K(4,4)& K(4,4)& K(4,4) },
\]
with $a$, $b$ and $e$ positive: $K(i,i) > K(4,4)$ for $i = 1,2,3$;
$K(i,i) > \sup\{a, b, e\}$ for $i = 1,2,3$; $\inf\{a,b,e\} > K(4,4)$.

\noindent For $a$, $ b$, $e$ distinct, $K$ is not symmetrizable and $K_{\vert _{\{
1,2,3\}\times\{1,2,3\}}}$ is its unique symmetrizable principal
submatrix of order $3$.
Moreover, the matrix $K$ is an inverse $M$-matrix.
\end{Remark}
%

\begin{Remark}
 For permanental vectors with symmetrizable kernel,
one might think that assuming the infinite divisibility of all its
triplets would lead to the infinite divisibility of the vector itself.
As it has been noticed in \cite{EK}, the Brownian sheet provides a
counter-example. Here is another one found in \cite{JS}. Indeed
the following matrix $B$ is a $4\times4$-covariance matrix of a
centered Gaussian vector $(\eta_1, \eta_2, \eta_3, \eta_4)$ such
that for every triplet of distinct indexes $i$, $j$, $k$, $(\eta_i^2,
\eta_j^2, \eta_k^2)$ is infinitely divisible, but $(\eta_1^2, \eta
_2^2, \eta_3^3, \eta_4^3)$ is not infinitely divisible:
\[
B = \pmatrix{ 1 & 0,50 & 0,35 & 0,40
\cr
0,50 &1& 0,50 & 0,26
\cr
0,35 & 0,50 &
1 & 0,50
\cr
0,40 & 0,26 & 0,50 & 1 }
\]
and $B^{-1}(2,4)$ is positive.
\end{Remark}

\begin{Remark} Let $(G(i,j))_{1 \leq i,j \leq n}$ be the kernel of
a permanental vector $\psi$. We assume that $G$ is nonsingular. For
$\alpha$ in $[0,1]$, consider now the $2n\times2n$-matrix $H(\alpha
)$ defined by
\[
H(\alpha) = \lleft[\matrix{ G & \alpha G
\cr
\alpha G & G } \rright] .
\]
The matrix $H(1)$ is the kernel of the vector $(\psi, \psi)$. The
matrix $H(0)$ is the kernel of the permanental vector $(\psi, \tilde
{\psi})$, where $\tilde{\psi}$ is an independent copy of $\psi$.
\end{Remark}

\begin{Proposition}
If $G$ does not contain any symmetrizable $3\times
3$-principal submatrix, then for any $\alpha$ in $(0,1)$, $H(\alpha
)$ is not the kernel of a permanental vector.
\end{Proposition}

\begin{pf} For $x$ complex number, we have
\begin{eqnarray*}
\operatorname{det} \bigl(H(\alpha) - xI \bigr) & = &\operatorname{det}
\lleft[\matrix{ G - xI & \alpha G
\cr
\alpha G & G - xI} \rright]
\\
& =& \bigl\llvert ( G - x I)^2 - \alpha^2
G^2 \bigr\rrvert = \bigl\llvert (1 + \alpha) G - xI \bigr\rrvert
\bigl\llvert (1- \alpha)G - xI \bigr\rrvert,
\end{eqnarray*}
since $\alpha G$ and $(G - xI)$ commute. Hence, $H(\alpha)$ satisfies
the first condition of Vere-Jones criterion of existence of a
permanental vector. Moreover, for $\alpha< 1$, $H(\alpha)$ is not singular.

By assumption, $G$ has no zero entry (if not it would contain a
symmetrizable $3\times3$-principal submatrix) and thanks to Theorem~\ref{T1}, it is hence diagonally equivalent to an inverse $M$-matrix.
In particular, there exists a signature matrix $\sigma$ such that the
entries of $\sigma G \sigma$ are all strictly positive. Consequently,
we can assume that the entries of $H(\alpha)$ are all strictly positive.

For $\alpha$ in $(0,1)$, $H(\alpha)$ is not an inverse $M$-matrix.
Indeed, write
\[
H(\alpha) = H = %
\lleft[\matrix{ H_{11} & H_{12}
\cr
H_{21} & H_{22} } \rright] %
,
\]
with $H_{11}= H_{22} = G$ and $H_{12}= H_{21} = \alpha G$, then
\[
H^{-1} = %
\lleft[\matrix{ (H/H_{22})^{-1}
& -(H/H_{22})^{-1} H_{12} (H_{22})^{-1}
\cr
-H_{22}^{-1} H_{21} (H/H_{22})^{-1}
& (H/H_{11})^{-1} } \rright], %
\]
where $H/H_{11}$ is the Schur complement of $H_{11}$ in $H$ defined by
\[
H/H_{11} = H_{22} - H_{21} H_{22}^{-1}
H_{12},
\]
and similarly $H/H_{22}$ is the Schur complement of $H_{22}$ in $H$:
\[
H/H_{22} = H_{11} - H_{12} H_{22}^{-1}
H_{21}.
\]
Then, as it has been noticed by Johnson and Smith \cite{JS}, $H$ is an
inverse $M$-matrix iff:
\begin{longlist}[(iii)]
\item[(i)] $H/H_{11}$ is an inverse $M$-matrix,

\item[(ii)] $H/H_{22}$ is an inverse $M$-matrix,

\item[(iii)] $(H_{22})^{-1} H_{21} (H/H_{22})^{-1} $ has nonnegative entries
only,

\item[(iv)] $(H/H_{22})^{-1} H_{12} (H_{22})^{-1}$ has nonnegative entries only.
\end{longlist}

For $H = H(\alpha)$ with $\alpha$ in $[0,1)$, this criterion gives
the following:

(i) and (ii): $(1 - \alpha^2) G$ is an inverse $M$-matrix.

(iii) and (iv): ${\alpha\over1 - \alpha^2} G^{-1}$ has only
nonnegative entries.

\noindent Hence unless $\alpha= 0$, $H(\alpha)$ is never an inverse $M$-matrix.

Now assume that there exists a permanental vector admitting
$H(\alpha)$ for kernel. Then we know that this permanental vector is
not infinitely divisible. But
as soon as $G$ does not contain any $3\times3$-covariance matrix,
neither does $H(\alpha)$.
Thanks to Theorem~\ref{T1}, a permanental vector that would admit
$H(\alpha)$ for kernel should be infinitely divisible. Hence,
$H(\alpha)$ can not be the kernel of a permanental vector.
\end{pf}

Note that for $G$ symmetric positive definite matrix, we know that
$H(\alpha)$ is still a covariance matrix. But the corresponding
$2n$-dimensional permanental vector is never infinitely divisible,
because Conditions (i) and (iii) above are always antagonistic.

%

\printaddresses
\end{document}